\documentclass[12pt,a4paper,twoside,final,notitlepage, reqno]{article}

\usepackage[english]{babel}
\usepackage[latin1]{inputenc}
\usepackage[a4paper,left=3.5cm,right=3.5cm]{geometry}
\usepackage{amssymb}
\usepackage{amsthm} 
\usepackage{amsmath}
\usepackage{mathtools}
\usepackage{amscd}
\usepackage{geometry}
\usepackage{graphics,graphicx}
\usepackage{epstopdf}
\usepackage[usenames, dvipsnames]{color}
\usepackage{hyperref}
\usepackage[textsize=small]{todonotes}
\usepackage{booktabs}
\usepackage{subfigure} 

\graphicspath{{figs/}}
\makeatletter
\def\input@path{{figs/}}
\makeatother

\theoremstyle{definition}
\newtheorem{thm}{Theorem}%[section]

\newtheorem{prop}{Proposition}

\theoremstyle{definition}
\theoremstyle{definition}
\newtheorem{rem}{Remark}%[section]
\newtheorem{exmp}{Example}%[section]

\providecommand{\keywords}[1]  {\textbf{Keywords:} #1}
\providecommand{\subjclass}[1] {\textbf{Subject classification:} #1}

\let\div\undefined
\DeclareMathOperator{\div}{div}

\DeclareMathOperator{\tr}{tr}

\usepackage{esdiff}
\renewcommand{\d}[1]{\, \mathrm{d} #1}

\newcommand{\R}{\mathbb{R}}

\renewcommand{\phi}{\varphi}

\addto\captionsfrench{}
\def\R{\mathbb{R}}

\def\nt{\nabla_\Gamma}
\def\na{\nabla}
\def\div{\mathrm{div}}

\def\Hk1{\mathrm{H}^{k+1}}
\def\HH1{\mathrm{H}^1}
\def\L2{\mathrm{L}^2}  
 
\def\LL{\mathrm{L}^}  
\def\Hexpo{\mathrm{H}^} 
\def\H2{\mathrm{H}^2}

\def\d{\mathrm{d}}
\def\nn{\boldsymbol{\mathrm{n}}}

\def\c{\mathcal{C}^} 
  
\def\P{\mathbb{P}^}

\def\omhh{\Omega_h}
\def\omh1{\omhh^{(1)}}
\def\omhr{\omhh^{(r)}}

\def\omgam{(\Omega,\Gamma)}

\def\ghh{\Gamma_h}
\def\gh1{\Gamma_h^{(1)}}
\def\ghr{\Gamma_h^{(r)}}

\def\ft{F_T}
\def\fte{\ft^{(e)}}
\def\ftre{F_{T^{(r)}}^{(e)}}
\def\ftr{\ft^{(r)}}
\def\ftdeux{\ft^{(2)}}

\def\dist{\mathrm{dist}}

\def\tref{\hat{T}}
\def\te{{T}^{(e)}}
\def\tdeux{{T}^{(2)}}
\def\tr{{T}^{(r)}}
\def\tauh{\mathcal{T}_h^{(1)}}
\def\taur{\mathcal{T}_h^{(r)}}
\def\taue{\mathcal{T}_h^{(e)}}

\usepackage{authblk}

\title{\bf{\Large{
      Numerical study of a diffusion equation with Ventcel boundary condition using curved meshes
    }}
}

\author[1]{
Fabien Caubet \thanks{fabien.caubet@univ-pau.fr}
}
\author[2]{
Joyce Ghantous \thanks{joyce.ghantous@univ-pau.fr}
}
\author[3]{
Charles Pierre \thanks{charles.pierre@univ-pau.fr}
}

\affil[1,2,3]{
  Laboratoire  de Math\'ematiques et de leurs Applications, 
  , UMR~CNRS~5142, \protect \\
  Universit\'e de Pau et des Pays de l'Adour, France.
}

\usepackage{fancyhdr}

\begin{document} 

\date{\today}

\maketitle
\noindent
\\ \\ 
\\ \\
% \acknow{}
% \\ \\
%
%
\begin{abstract}
In this work is provided a numerical study of a diffusion problem involving a second order term on the domain boundary (the Laplace-Beltrami operator) referred to as the \textit{Ventcel problem}.
A variational formulation of the Ventcel problem is studied, leading to a finite element discretization.
The focus is on the resort to high order curved meshes for the discretization of the physical domain.
The computational errors are investigated both in terms of geometrical error and of finite element approximation error, respectively associated to the mesh degree $r\ge 1$ and to the finite element degree $k\ge 1$. 
The numerical experiments we led allow us to formulate a conjecture on the \textit{a priori} error estimates depending on the two parameters $r$ and $k$. 
In addition, these error estimates rely on the definition of a functional \textit{lift}  with adapted properties on the boundary to move numerical solutions defined on the computational domain to the physical one. 
\end{abstract}

\section*{Introduction}

\paragraph{Motivation.} 
On the one hand, in various situations, we have to numerically solve a problem (a partial differential equation) on a non-polygonal geometry. This~requires~the~use~of~high order meshes in order to well approximate it. On the other hand, in several industrial applications, objects or materials surrounded by a thin layer with potentially other properties (typically a surface treatment or corrosion) have to be considered. The presence~of~this~layer~causes some difficulties while discretizing the domain and numerically solving the problem. To overcome this problem, the domain is approximated asymptotically by an other one without a thin layer but equipped with artificial boundary conditions, like {\it Ventcel boundary condition}. The physical properties of the thin layer are then contained in the boundary condition.
\vspace{0.1cm}

%{\color{red}
This paper focuses on the resolution of a problem involving higher order boundary condition and numerically evaluates the \textit{a priori} error produced by a finite element approximation on higher order meshes, distinguishing the {\it geometrical error} from the {\it approximation error}.%}
%The main focus of this paper is to consider a system that is equipped with this higher order boundary condition, and after that assess the \textit{a priori} error produced by a finite element approximation, on higher order meshes.

\paragraph{The Ventcel problem and its approximation.} 
Let  $\Omega$ be a domain in $\R^{d}$, $d=2$, 3, 
 with a smooth boundary $\Gamma$. 
The Laplace-Beltrami operator on $\Gamma$ is denoted by $\Delta_\Gamma$.
Relatively to the source terms $f$ and $g$ and to the constants $\kappa \ge 0$, $\alpha,\,\beta>0$, the Ventcel problem reads,
\begin{equation}
\label{1}
  \left\{
      \begin{aligned}
       -\Delta u + \kappa u &=& f &  \quad  \text{ in } \, \Omega ,\\
        -\beta \Delta_{\Gamma} u + \partial_{\nn} u + \alpha u &=& g &  
        \quad   \text{ on } \, \Gamma ,\\
      \end{aligned}
    \right.
\end{equation}
with $\nn$ the external unit normal to $\Gamma$ 
and $\partial_{\nn} u$ the normal derivative of $u$ along $\Gamma$. 
The theoretical properties of the solution of  problem 
\eqref{1} 
have been studied in \cite{ventcel1}.
% This problem admits a unique solution $u \in \HH1 \omgam:= \{ u \in \Hexpo{1}(\Omega), \ u_{|_\Gamma} \in \Hexpo{1}(\Gamma) \}$ (see theorem \ref{th_existance_unicite_u}). We here focus on an efficient numerical approximation of the solution of such a problem.%\vspace{0.1cm}

Due to the presence of the second order  
%Laplace-Beltrami
term in the boundary condition,
the domain $\Omega$ is required to be smooth and thus non-polygonal:
from the numerical point of view, the computational domain $\Omega_h$ (the mesh domain) will not fit the physical one: $\Omega_h \neq \Omega$.
% Due to the presence of a higher order operator in the boundary condition, 
% more accuracy is required when discretizing the domain $\Omega$. Hence, we construct higher order meshes of geometrical order $r \ge 1$ by the help of a transformation that maps the reference simplex onto a curved element (see section \ref{mesh}). 
In a context of finite element methods 
of high order  $k \ge 2$,
it then is necessary to resort to high order meshes of 
geometrical degree $r \ge 2$ to preserve the numerical solution's accuracy. Some methods have been widely studied, see, \textit{e.g.},
% was developed and modified by many authors
\cite{elliott,ed,dubois,scott,nedelec}. 
% Additionally, we use a $\P k$ finite element method with a degree $k \ge 1$ to solve the considered problem. 
% Notice that the isoparametric approach, which is the case of $r=k$, is widely treated in similar cases (see \cite{elliott,ed}).%\vspace{0.1cm}

The approximation of the Laplace equation on a surface 
has been studied in this framework 
by Demlow \textit{et al.} in \cite{D1,D2}. 
In these works, a distinction is made  between the geometrical error 
induced by the setting of the computational domain $\Omega_h \neq \Omega$
and the approximation error related to the finite element method.
The purpose of this approach is to highlight the influence of 
the geometrical degree $r$ of the mesh and the finite element approximation degree $k$  on the total computational error. 
Thereby, one can assess which is the optimal degree of the finite element method $k$ to chose depending on the choice of the geometrical degree $r$. 
%Note that the case $r=1$ is treated in \cite{D2}.%\vspace{0.1cm}

In the present context where $\Omega_h \neq \Omega$, a crucial issue arises: how does one compare the numerical solutions $u_h$ to the exact one, in order to derive \textit{a priori} error estimates? 
To circumvent this, a lift of  $u_h$ onto $\Omega$ is defined: 
in \cite{dubois}, Dubois introduced such a lift based 
on the orthogonal projection  onto the boundary $\Gamma$, which 
further was improved in terms of regularity  by Elliott \textit{et al.} \cite{elliott}. This lift however does not fit the orthogonal projection on the computational domain boundary.
An alternative definition is introduced in this paper which will be used to perform a numerical study of the computational error of problem \eqref{1}.

\paragraph{Paper organization.} 
In section $\ref{sec22}$, 
after introducing some general mathematical tools, 
is stated and proven the well-posedness of the Ventcel problem \eqref{1}. 
The following section $\ref{mesh}$ is devoted to the definition of the curved meshes of $\Omega$. 
In section \ref{sec:num-exp} are presented the discretization 
of the Ventcel problem \eqref{1}, the lift operator which is the keystone of the \textit{a priori} error estimations and numerical experiments studying the method convergence rate depending on the mesh geometrical degree $r$ and on the finite element approximation degree $k$.
The paper ends with a conclusion section presenting our conjecture on \textit{a priori} error estimates.

\section{Study of the Ventcel problem}
\label{sec22}
\paragraph{Some mathematical tools.}
Let us 
denote $\Omega$  a bounded connected open subset of $\R^{d}$ %$(d=2,3)$, 
with a smooth boundary $\Gamma:=\partial{\Omega}$ 
 at least of $\c2$ regularity.
The  unit normal to $\Gamma$ pointing outwards is denoted by $\nn$.
%and $\partial_{\nn} u$ is the normal derivative along $\Gamma$ 
%of a function $u$ on $\Omega$.
The classical spaces $\LL 2(\Omega)$, $\LL 2(\Gamma)$,  $\Hexpo{1}(\Omega)$ and $\Hexpo{1}(\Gamma)$ are considered 
and we introduce 
the following Hilbert space and its associated norm (see \cite[Lemma 2.5]{ventcel1})
\begin{equation*}
    \Hexpo{1}\omgam := \{ u \in \Hexpo{1}(\Omega), \ u_{|_\Gamma} \in \Hexpo{1}(\Gamma) \},\qquad 
    \|u\|^2_{\Hexpo{1}\omgam} :=  \|u\|^2_{\Hexpo{1}(\Omega)} 
+ \|
u_{|_\Gamma}\|^2_{\Hexpo{1}(\Gamma)}.
\end{equation*}
We consider the classical surface operators 
(see, e.g., \cite[p. 192-196]{livreopt}): 
\begin{itemize}
\item the {\it tangential gradient} of 
  $u\in\HH1(\Gamma)$ is given by $\nt u :=\na \tilde{u}  - (\na \tilde{w} \cdot \nn)\nn$, where $\tilde{u} \in \HH1(\R^d)$ is any extension of $u$;
\item the {\it tangential divergence} of $W\in\HH1(\Gamma,\R^d)$ is  
  $\div_{\Gamma} W : = 
  \div \tilde{W}- (\mathrm{D}\tilde{W}\, \nn)\cdot \nn$, 
  where~$\tilde{W} \in \HH1(\R^d, \mathbb{R}^d)$ 
  is any extension of $W$ and $\mathrm{D}\tilde{W}$ is the differential of $\tilde{W}$;
\item the {\it Laplace-Beltrami operator} of 
  $u\in \Hexpo{2}(\Gamma)$ is given by $\Delta_\Gamma u  := \div_\Gamma (\nt u)$.
\end{itemize} 
Finally, the 
following fundamental result is recalled, see, \textit{e.g.}, \cite{tubneig} and \cite[\S 14.6]{GT98}.

\begin{prop}
\label{tub-th}
Let $\Omega$ and  $\Gamma= \partial \Omega$ be as stated previously. Let $\d : \R^d \to \R$ be the signed distance function with respect to $\Gamma$ defined by,
\begin{displaymath}
  \d(x) :=
  \left \{
    \begin{array}{ll}
      -\dist(x, \Gamma)&  {\rm if } \, x \in \Omega ,
      \\
      0&  {\rm if } \, x \in \Gamma ,
      \\
      \dist(x, \Gamma)&  {\rm otherwise},
    \end{array}
  \right. \qquad 
  {\rm with} \quad \dist(x, \Gamma) := \inf \{|x-y|,~ \ y \in \Gamma \}.
\end{displaymath}
Then there exists a tubular neighborhood $\mathcal{U}_{\Gamma} $ 
of $\Gamma$ where {$\d$ is a $\c2$ function}. Its gradient 
$\na \d$ %$\nn:=\na \d$ 
is an extension of the external unit normal $\nn$ to $\Gamma$. Additionally, in this neighborhood $\mathcal{U}_{\Gamma}$, the orthogonal projection~$b$ onto $\Gamma$ is uniquely defined
and given by
\begin{displaymath} 
b\, :~ x \in \mathcal{U}_{\Gamma}  
\longmapsto    b(x):=x-\d(x)
\na \d(x) \in \Gamma.
\end{displaymath}
\end{prop}

\paragraph{Well-posedness of problem \eqref{1}.}
The weak form of \eqref{1} is to find $u \in \HH1 \omgam$ such that,
\begin{equation}
  \label{fv1}
  \forall \ v \in \HH1 \omgam, \quad 
    a(u,v)= l(v) := \int_{\Omega} f v \d x +\int_{\Gamma} g v \d \sigma,
\end{equation}
where $a(\cdot,\cdot)$ is defined on $\HH1\omgam^2$ by:
\begin{equation}
  \label{eq:def-a}
  a(u,v) := \int_{\Omega} \nabla u \cdot \nabla v \d x +\kappa \int_{\Omega}  u  v \d x + \beta \int_{\Gamma} \nabla_{\Gamma} u \cdot \nabla_{\Gamma} v \d \sigma + \alpha \int_{\Gamma} u  v \d \sigma.
\end{equation} 
Notice that the weak form \eqref{fv1} is equivalent to the system introduced in \eqref{1} as it was proven in \cite{ventcel1}.
\begin{thm}
\label{th_existance_unicite_u}
Let $\Omega$ and  $\Gamma= \partial \Omega$ be as stated previously. 
Let $\alpha$, $\beta >0$, $\kappa \ge 0 $, and $f \in \L2(\Omega)$, $g \in \L2(\Gamma)$. Then there exists a unique solution 
$u \in \Hexpo{1}\omgam$ to problem (\ref{fv1}).
\end{thm}
%%
%%
%%
%%
%%
%%
%In addition, it is proven in \cite[th. 3.2]{ventcel1} without getting into details that there exists a unique solution for the Ventcel problem. After that a quick proof is presented stating that a weak solution to (\ref{fv1}) is a strong solution of (\ref{1}). In \cite[th. 3.3]{ventcel1}, it is proven that
%there exists a (source term independent) constant $c>0$ such that, 
%\begin{displaymath}
%    \|u\|_{\H2\omgam} \le c ( \|f\|_{\L2(\Omega)} + \|g\|_{\L2(\Gamma)}).
%\end{displaymath}
%%%
%%%
%%%
%\\
%%%
%%%
%%%
%%%
%%%
%%%
%The proof of theorem \ref{th_existance_unicite_u}, though classical,
%is detailed here for the sake of completeness.
The proof of this theorem is classical and is briefly given in~\cite[th. 3.2]{ventcel1}. We detail it here for the sake of completeness. Let us notice that, additionally, it is proven in~\cite[th. 3.3]{ventcel1} that
there exists a (source term independent) constant $c>0$ such that
\begin{displaymath}
    \|u\|_{\H2\omgam} \le c ( \|f\|_{\L2(\Omega)} + \|g\|_{\L2(\Gamma)}).
\end{displaymath}
\begin{proof}
The proof relies on the Lax-Milgram theorem.
The linear form $l(\cdot )$ in (\ref{fv1}) and the bilinear form $a(\cdot ,\cdot )$ in (\ref{eq:def-a}) being continuous respectively on $\Hexpo{1}\omgam$ and on $\Hexpo{1}\omgam^2$, 
it remains to show that $a$ is coercive. 
We must distinguish between two cases.

%\textbf{Case 1.} If $\kappa \ne 0$, the result is obvious: 
\textbf{1) If $\kappa \ne 0$.} The result is obvious: 
for all $u\in \HH1 \omgam$, $a(u,u) \ge \min\{1, \kappa, \alpha, \beta \} \| u \|_{\HH1 \omgam}^2$.

%\textbf{Case 2.} If $\kappa = 0 $, we proceed by contradiction 
\textbf{2) If $\kappa = 0 $.} We proceed by contradiction assuming that there exists a sequence $(u_n)_{n\in \mathbb{N}^*}$ in $ \HH1 \omgam$ such that for all $n\ge 1$,
\begin{displaymath}
%\label{3}
     \|\nabla u_n\|^2_{\L2(\Omega)} + \beta \|\nabla_{\Gamma} u_n\|^2_{\L2(\Gamma)}+ \alpha \|u_n\|^2_{\L2(\Gamma)} < \frac{1}{n} \big( \|u_n\|^2_{\HH1(\Omega)} + \|u_n\|^2_{\HH1(\Gamma)} \big).
\end{displaymath}
It follows that $u_n \ne 0$ for all $n \ge 1$. Thus $u_n$ can be renormalized such that
%\begin{equation*}
%\|u_n\|_{\HH1\omgam}=1 \ \ \ \ \mbox{ and } \ \ \ \ \ 
%    \|\nabla u_n\|^2_{\L2(\Omega)} + \beta \|\nabla_{\Gamma} u_n\|^2_{\L2(\Gamma)}+ \alpha  \|u_n\|^2_{\L2(\Gamma)} < \frac{1}{n} . %, \ \ \ \ \ \forall \ n\ge 1.
%\end{equation*}
$\|u_n\|_{\HH1\omgam}=1$ and it satisfies $\|\nabla u_n\|^2_{\L2(\Omega)} + \beta \|\nabla_{\Gamma} u_n\|^2_{\L2(\Gamma)}+ \alpha  \|u_n\|^2_{\L2(\Gamma)} < \frac{1}{n} $.
Therefore
%\begin{align}
%    \label{cv1}
%    \nabla u_n \to 0  \ \ \ \mbox{ in } \ \L2(\Omega), \ \ \ \ \ 
%\nabla_\Gamma u_n \to 0  \ \ \ \mbox{ in } \ \L2(\Gamma), 
%\end{align}
%\begin{equation}
%    \label{cv3}
%u_n \to 0 \ \ \ \mbox{ in } \ \L2(\Gamma).
%\end{equation}
\begin{equation} 
\label{cv1et2}
    \nabla u_n \to 0  \ \ \ \mbox{ in } \ \L2(\Omega), \quad  \nabla_\Gamma u_n \to 0  \ \ \ \mbox{ in } \ \L2(\Gamma) \quad\mbox{ and } \quad  u_n \to 0 \ \ \ \mbox{ in } \ \L2(\Gamma).
\end{equation}
Since $(u_n)_n$ is bounded in $\HH1\omgam$, there exists $u\in \HH1 \omgam $ such that
%\begin{equation*}
%\label{cv4}
%    u_n  \rightharpoonup u  \ \ \mbox{in} \ \ \HH1\omgam,
%\end{equation*}
$u_n  \rightharpoonup u $ in~$ \HH1\omgam$,
and since ${\HH1\omgam \hookrightarrow~\L2\omgam}$ is a compact injection, we obtain
\begin{equation}
\label{cv5}
    u_n \to u \ \ \mbox{in} \ \ \L2\omgam.
\end{equation}
%For all $n\ge 1$, since $\|u_n\|_{\HH1\omgam}=1$, we have
%\begin{equation*}
%     \|\nabla u_n\|^2_{\L2(\Omega)} + \|\nabla_{\Gamma} u_n\|^2_{\L2(\Gamma)}+  \|u_n\|^2_{\L2(\Gamma)} + \|u_n\|^2_{\L2(\Omega)}=1.
%\end{equation*}
%Then, by passing to the limit and using the convergences given in \eqref{cv1et2} and \eqref{cv5}, we obtain 
Passing to the limit in $\|u_n\|^2_{\HH1\omgam} = \|\nabla u_n\|^2_{\L2(\Omega)} + \|\nabla_{\Gamma} u_n\|^2_{\L2(\Gamma)}+  \|u_n\|^2_{\L2(\Gamma)} + \|u_n\|^2_{\L2(\Omega)}=1$, 
and using the convergences given in \eqref{cv1et2} and \eqref{cv5}, we obtain 
%\begin{equation*}
$   \|u\|^2_{ \L2\omgam} = 1$.
%\end{equation*}
%
%The weak convergence $u_n  \rightharpoonup u \mbox{ in } \HH1\omgam$ implies that $ \nabla u_n \rightharpoonup \na u$ in~$\L2(\Omega).$
However, since $ \nabla u_n \rightharpoonup \na u$ in~$\L2(\Omega)$, we use~\eqref{cv1et2} and the uniqueness of the limit % in~$\L2(\Omega)$, we have 
to obtain $\na u = 0$ and, since~$\Omega$ is a connected set, it follows that $u=C \in \mathbb{R}$. 
Finally, $u_n \to u$ in~$\L2(\Gamma)$ and also $ u_n \to 0$ in~$\L2(\Gamma)$, these two points yield  $u = 0 = C$ which contradicts  $\|u\|_{\L2\omgam}=1$ and concludes the proof of the coercivity. 
\end{proof}

\section{Curved mesh definition}
\label{mesh}
In this section are defined curved meshes of geometrical degree  $r\ge 1$ of the domain~$\Omega$.
From now on, the domain $\Omega\subset\R^d$, $d=$2 or 3, is assumed to be at least $\c {r+2}$ regular, and~$\tref$ denotes the reference simplex of dimension $d$.
The definition steps are the following (see~\cite{elliott,scott,dubois} for more details).
\begin{enumerate}
    \item Construct an affine mesh $\tauh$ of $\Omega$ composed of simplexes  $T$. 
    \item For each $T \in \tauh$, a mapping 
      $\fte:~\tref\rightarrow \te := \fte(\tref)$ is designed, 
      so that the {\it exact element} $\te$ form a curved mesh $\taue$ whose domain exactly fits $\Omega$.
    \item For each $T \in \tauh$, the mapping  $\fte$ is interpolated by a polynomial $\ftr$ of degree $r$.
The associated elements  $\tr := \ftr(\tref)$ form a curved mesh $\taur$
of degree $r$ of $\Omega$.
\end{enumerate}
%This method and elements are detailed in  \cite{elliott,scott,dubois}.
%%
%%
%%
%%
%%
%%
%%
%%
%%
%%
%%
%%
\paragraph{Affine mesh.}
Let $\tauh$ be a %polyhedral 
mesh of $\Omega$ made of simplexes of dimension $d$ (triangles or tetrahedra), 
it is chosen as quasi-uniform and henceforth shape-regular (see \cite[def. 4.4.13]{quasi-unif}). % and~\cite[p. 5]{D4}). 
The mesh domain is denoted by $\omh1:= \cup_{T\in  \tauh}T$ and its boundary by $\gh1 :=\partial \omh1$, which is composed of $(d-1)$-dimensional simplexes that form a  mesh of $\Gamma = \partial \Omega$. The vertices of $\gh1$ are assumed to lie on $\Gamma$. We define the mesh size $h:= \max\{{\rm diam}(T); T \in \tauh \}$. 
%\\
To each  $T \in \tauh$ is associated an affine function 
 $\ft : \tref \to T = \ft(\tref)$. 

\begin{rem}%{As in \cite{elliott}, 
For a sufficiently small $h$, the mesh boundary satisfies 
$\gh1 \subset \mathcal{U}_\Gamma$, 
where $\mathcal{U}_{\Gamma}$ is the tubular neighborhood 
given in proposition~\ref{tub-th}.
This guaranties that the orthogonal projection~$b: \gh1\rightarrow \Gamma$ is one to one which is 
required for the construction of the exact mesh.
 \end{rem}
\begin{exmp}
%\label{expl}
%\emph{
In the two dimensional case is displayed  the case of a triangle $T \in \tauh$, with $T \cap \Gamma = \{ v_1,v_2\}$, together  with the mapping $\ft$ that maps  $\tref$ into $T$.
\begin{center}
\vspace{-0.05cm}
\begin{tikzpicture}[scale=0.81]
%%%% conctruction of the reference element

\draw (0.5,0.5) node[above] {$\hat{T}$};

\draw (0,0) node  {$\bullet$};
\draw (2,0) node  {$\bullet$};
\draw (0,2) node  {$\bullet$};

\draw (0,0) node[below] {$\hat{v}_1$};
\draw (2,0) node[below]  {$\hat{v}_2$};
\draw (0,2) node[above]  {$\hat{v}_3$};

\draw (0,0) -- (2,0) ;
\draw (0,0) -- (0,2);
\draw (2,0) -- (0,2) ;

%%%%%%%%%%%%%%%%%%% curved ref element 

%\draw[red] (1,0) node  {$\bullet$};
%\draw[red] (1,1) node  {$\bullet$};
%\draw[red] (0,1) node  {$\bullet$};

%\draw (1,0) node[below] {$e_3$};
%\draw (1,1) node[above]  {$e_1$};
%\draw (0,1) node[left]  {$e_2$};

%%%%%%%%%%%%%%%%%%%%%%%%%%%%%%%%%%%%% arrow and function

\draw[thick, blue] [->] (2.8,1) -- (4.3,1);
\draw[blue] (3.5,1) node[above] {${\ft}$};

%%%%%%%%%%%%%%% curved case : transformation
%\draw[red] (3.5,1) node[above] {$\textcolor{red}{\ftr}$};

%%%%%%%%%%%%%%%%%%%% GAMMA %%%%%%%%5

%\draw (6,0) arc (180:90:3);
%\draw (9,3.2) node {$\Gamma$};

%%%%%%%%%%%%%%%%%%%% triangle construction 

%%%%%%%%%%%%%%%%%%%%%%%%%%%%%%%%%%%%% name  

%%%%%%%%%%% linear
\draw (7.5,0.8) node[above] {$T$};

%%%%%%%%%%% curved 
%\draw[red] (7.4,0.8) node[above] {$\tr$};

%%%%% points %%%%%%%%%%%%%%%%%%%%%

\draw (6,0) node  {$\bullet$};
\draw (7.36,2.5) node  {$\bullet$};
\draw (9,0) node  {$\bullet$};

\draw (6,0)  node[left] {$v_1$};
\draw (7.36,2.6) node[above]  {$v_2$};
\draw (9,0) node[below]  {$v_3$};

\draw (7.36,2.57) -- (9,0);
\draw (9,0) -- (6,0);

%%%%%%%%%%%%%%%%%%%%%%%%%%% curved triangle

%\draw (6,0)  node[below] {$v\r1$};
%\draw (7.36,2.6) node[above]  {$v\r3$};
%\draw (9,0) node[below]  {$v\r2$};
%\draw (6.3,1.3) node[left]  {$e\r2$};
%\draw[red] (6.3,1.3) node  {$\bullet$};
%\draw (7.5,0) node[below]  {$e\r3$};
%\draw (8.15,1.3) node[right]  {$e\r1$};
%\draw[red] (7.5,0) node  {$\bullet$};
%\draw[red] (8.15,1.3) node  {$\bullet$};

%%%%%%%%%%%%%%%%%%%%%%%%% the linear triangle

\draw[blue] (6,0) -- (7.35,2.5);

%%%%%%%%%%%% x in T_1
\draw (6.7,1.3) node  {$\bullet$};
\draw (6.8,1.3) node[right]  {$x$};

%%%%%%%%%%% b(x)
\draw (6.4,1.4) node  {$\bullet$};
\draw (6.3,1.35) node[left]  {$b(x)$};

%%%%%%%%%% linking x and b(x)

\draw[purple] (6.9,1.25) -- (6,1.55);

%%%%%%%%% angle droit

\draw[purple] (6.475,1.2) -- (6.53,1.38);
\draw[purple] (6.475,1.2) -- (6.35,1.238);

%%%%%%%%%%%%%%%%%%%%%%% the tubular neighborhood

%\draw[blue] (6.4,0) arc (180:90:2.5);

%\draw[blue] (5.6,0) arc (180:90:3.2);

%%%%%%%%%%%%%%%%%%%%%% GAMMA NOT EXACT

\draw plot [domain=-0.2:1.5] (\x+6,-1.786*\x^2+4.3*\x);

\draw (5,-0.6) node {$\Gamma$};
\draw[thick] [->] (5.2,-0.6) -- (5.88,-0.5);

%%%%%%%%%%%%%%%%%%%%%% GAMMA_h

\draw[blue] (5,2.6) node {$\gh1$};
\draw[thick, blue] [->] (5.2,2.4) -- (7,1.8);

\end{tikzpicture}
\end{center}
%\emph{
% Notice that $[v_1,v_2] = T\cap\gh1 \subset \mathcal{U}$ 
% and that the projection $b$ is one to one from $[v_1,v_2]$ to of any point $x \in [v_1,v_2]$ is uniquely defined on the border $\Gamma$, by Proposition~\ref{tub-th}.
\end{exmp}
\paragraph{Exact mesh $\taue$.}
After the early works of Scot \cite{scott} 
and  Lenoir \cite{Lenoir1986} defining 
transformations towards curved elements, Dubois \cite{dubois} first introduced a definition based on the orthogonal projection $b$ onto $\Gamma$, 
further developed by Elliott \textit{et al.} \cite[\S 4]{elliott} 
in terms of regularity, which definition is recalled here.

Let us first point out that, because of the quasi uniform assumption 
made on the mesh $\tauh$, and for $h$ sufficiently small,
a mesh element $T$ cannot have $d+1$ vertices on the boundary $\Gamma$.
We define internal elements as those having at most one  vertex on the 
boundary $\Gamma$, whereas other elements have:
\begin{itemize}
\item 2 vertices on the boundary in the two dimensional case;
\item 2 or 3 vertices on $\Gamma$ in the 3D case,
  forming either an edge or a face respectively.
\end{itemize}
The case of internal elements is skipped by setting $\fte=\ft$.

Let then $T\in\tauh$ a non-internal element, denote $v_i = \ft(\hat{v}_i)$ its vertices, $\hat{v}_i$ being the vertices of $\tref$, and define $\varepsilon_i=1$  if $v_i\in  \Gamma$ or  $\varepsilon_i=0$ otherwise.
To $\hat{x}\in \tref$ is associated its barycentric coordinates $\lambda_i$ associated to the vertices $\hat{v}_i$ of $\tref$ .
We introduce 
%relatively to $\hat{x}\in \tref$,
$  \lambda^\star := \sum_{i=1}^{d+1} \varepsilon_i \lambda_i $
and
$  \hat{y} := \dfrac{1}{\lambda^\star}\sum_{i=1}^{d+1} \varepsilon_i \hat{v}_i\in\tref$.
The mapping $\fte:~\tref \rightarrow \te$ is given by,
\begin{equation}
  \label{eq:def-fte}
  \fte(\hat{x}) := x + (\lambda^\star)^{r+2} ( b(y) - y),
  \quad  {\rm with} \quad 
  x = \ft(\hat{x})
  \quad  {\rm and} \quad 
  y = \ft(\hat{y}).
\end{equation}
\begin{rem}
\label{rem:Fe_T}
For $x\in T \cap \ghh$, we have that $\lambda^\star = 1$ 
and so $y = x$
inducing that
$\fte(\hat{x}) = b(x)$: 
$\fte \circ \ft^{-1} = b$ on $T \cap \ghh$ which is then mapped on $\Gamma$ following the orthogonal projection $b$.
%\\
The mapping $\fte$ has been shown in \cite{elliott} to be $C^{r+1}$ regular on $\tref$.
\end{rem}

% For any $T \in \tauh$, we define the following $\c1$ transformation $\fte$ that maps the reference element $\tref$ to an exact element $\te$  
% \begin{equation*}
% %\fonction{\fte}{\tref  }{\te :=\fte( \tref) }{ \hat{x}}{\displaystyle  \fte(\hat{x}) =\ft (\hat{x})+\rho_T(\hat{x});}
% \fte \, : \, \hat{x} \in \tref  \longmapsto  \ft (\hat{x})+\rho_T(\hat{x}) \in \te :=\fte( \tref)  ,
% \end{equation*}
% where $\ft: \tref \to T$ is given above and $\rho_T: \tref \to \te$ is a given $\c1$ mapping. The role of $\rho_T$ is to move or translate the points of $T$ to form the exact triangle $\te \subset \Omega$. In \cite[\S 4]{elliott} and \cite[\S 3.2]{ed}, the construction of $\rho_T$ is given step by step. Scott \cite{scott} gives an explicit construction of an exact triangulation in two dimensions which is generalized by Lenoir \cite{Lenoir1986}. Here, in this paper, we got inspired by Dubois's work \cite{dubois} which uses the normal projection \eqref{db} on a 3d problem. Thus, this transformation $\fte$ combines both the projection b and the definition in~\cite{Lenoir1986,scott}. It has to be noted that $\fte=\ft$ on any internal element $T$, by definition of $\fte$. 

\begin{exmp}
%\emph{
Consider three triangles $T_1$, $T_2$ and $T_3$ in $\R^2$ 
%with vertices $(v_i)_{i=1}^5$ 
as displayed below. For $i=1,2,3$, we have the following transformation $F_{T_i}^{(e)}\circ F_{T_i}^{-1}$ that maps $T_i \in \tauh$ into $\te_i$ as follows,
% }
\begin{center}
\vspace{-0.1cm}
\begin{tikzpicture}[scale=0.90]
%%%%%%%%%%%%%%%%%%%%%% GAMMA NOT EXACT

\draw plot [domain=-0.2:1.5] (\x-1,-2*\x^2+4.4*\x);

\draw (-2.3,-0.6) node {$\Gamma$};
\draw[thick] [->] (-2,-0.6) -- (-1.12,-0.5);

%%%%%%%%%% GAMMA %%%%%%%%%%%%%

%\draw (2.3,3) node {$\Gamma$};
%\draw (-1,0) arc (180:90:3);

%%%%%%%%%%%%% T1 %%%%%%%%%%%%%%%
\draw (0,0) node  {$\bullet$};
\draw (0.18,2.4) node  {$\bullet$};
\draw (0.75,1.25) node  {$\bullet$};
\draw  (0,0) -- (0.75,1.25) ;
\draw  (0,0) -- (0.18,2.4);
\draw (0.75,1.25) -- (0.18,2.4) ;

\draw (0.4,1) node[above] {$T_1$};

%%%%%%%%%%% T2 %%%%%%%%%%%%%%%%%%%
\draw (0,0) node  {$\bullet$};
\draw (2,0) node  {$\bullet$};
\draw (0.75,1.25) node  {$\bullet$};

\draw (0,0) -- (2,0);
\draw (2,0) -- (0.75,1.25);
\draw (0.75,1.25) -- (0,0);

\draw (0.9,0.3) node[above] {$T_2$};

%%%%%%%%%%% T3 %%%%%%%%%%%%%%%%%%

\draw (-1,0) node  {$\bullet$};
\draw (0.18,2.4) node  {$\bullet$};

\draw[magenta] (-1,0) -- (0.18,2.4);
\draw (0,0) -- (-1,0);

\draw (-0.2,0.4) node[above] {$T_3$};

%%%%%%%%%% border %%%%%%%%%%%%%

%\draw[magenta] (-2,2.45) node {$\partial T_3$};
%\draw[magenta, thick] [->] (-1.8,2.2) -- (-0.2,1.62);

%%%%%%%%%%%%%% vertices%%%%%%%%%%%%%%

\draw (0.18,2.4) node[above] {$v_1$};
\draw (-1,0) node[left] {$v_2$};
\draw (2,0) node[below] {$v_3$};

\draw (0,0) node[below] {$v_4$};
\draw (0.8,1.25) node[right] {$v_5$};

%%%%%%%%%%%%%%%%%% transformation %%%%%%%%%%%%%

\draw[thick, magenta] [->] (3.2,1) -- (4.8,1);
\draw[magenta] (4,1) node[above] {$\textcolor{magenta}{F_{T_i}^{(e)}\circ F_{T_i}^{-1}}$};

%%%%%%%%%%%%%%%%%% second triangle %%%%%%%%%%%%
%%%%%%%%%%%%%%%%%%%%%% GAMMA NOT EXACT

\draw plot [domain=-0.2:1.5] (\x+6,-2*\x^2+4.4*\x);

\draw (4.7,-0.6) node {$\Gamma$};
\draw[thick] [->] (5,-0.6) -- (5.88,-0.5);

%%%%%%%%%% GAMMA_h %%%%%%%%%%%%%

\draw[magenta] plot [domain=-0:1.2] (\x+6,-2*\x^2+4.4*\x);

%\draw[magenta] (5,2.45) node {$\partial T^{(e)}_3$};
%\draw[magenta, thick] [->] (5.2,2.2) -- (6.48,1.62);

%%%%%%%%%% GAMMA %%%%%%%%%%%%%

%\draw (9.3,3) node {$\Gamma$};
%\draw (6,0) arc (180:90:3);

%%%%%%%%%%%%% T1 %%%%%%%%%%%%%%%
\draw (7,0) node  {$\bullet$};
\draw (7.18,2.4) node  {$\bullet$};
\draw (7.75,1.25) node  {$\bullet$};

\draw (7,0) -- (7.18,2.4);
\draw (7.18,2.4) -- (7.75,1.25);
\draw (7.75,1.25) -- (7,0);

\draw (7.4,0.8) node[above] {$T^{(e)}_1$};

%%%%%%%%%%% T2 %%%%%%%%%%%%%%%%%%%
\draw (7,0) node  {$\bullet$};
\draw (9,0) node  {$\bullet$};
\draw (7.75,1.25) node  {$\bullet$};

\draw (7,0) -- (9,0);
\draw (9,0) -- (7.75,1.25);
\draw (7.75,1.25) -- (7,0);

\draw (7.9,0.2) node[above] {$T^{(e)}_2$};

%%%%%%%%%%% T3 %%%%%%%%%%%%%%%%%%

\draw (6,0) node  {$\bullet$};
\draw (7.18,2.4) node  {$\bullet$};

%\draw (6,0) -- (7.18,2.4);
\draw (7,0) -- (6,0);

\draw (6.65,0.4) node[above] {$T^{(e)}_3$};

%%%%%%%%%%%%%% vertices%%%%%%%%%%%%%%

\draw (7.18,2.4) node[above] {$v_1$};
\draw (6,0) node[left] {$v_2$};
\draw (9,0) node[below] {$v_3$};

\draw (7,0) node[below] {$v_4$};
\draw (7.8,1.3) node[right] {$v_5$};
\end{tikzpicture}
\end{center}
%\emph{
$T_1$ and $T_2$ are internal and so are unchanged whereas $T_3$ 
(having 2 vertices on $\Gamma$) 
is not internal and mapped into a curved triangle 
with an edge exactly fitting $\Gamma$.
% $T_1$ has three vertices of which only $v_1$ is on the boundary $\Gamma$. Then $F_{T_1}^{(e)}\circ F_{T_1}^{-1}=Id$, since $T_1$ is an internal simplex of the mesh. Similarly, $F_{T_2}^{(e)}\circ F_{T_2}^{-1}=Id$ since $T_2 \cap \Gamma = \emptyset $. On the other hand, $T_3$ has two vertices $v_1, v_2$ on $\Gamma$ then $F_{T_3}^{(e)}\circ F_{T_3}^{-1} \ne Id$, and $\fte\circ \ft^{-1}([v_1,v_2]) \subset \Gamma$.
%}
\end{exmp}
%We define now the higher-order polynomial meshes $\omhh^{(r)}$ while using the previous transformation. 
%%
%%
%%
%%
%%
%%
%%
%%
\paragraph{Curved mesh $\taur$ with degree $r$.}
Let $T \in~\tauh$ and $r\ge 1$,
the exact mapping $\fte$ in 
\eqref{eq:def-fte} 
is interpolated as a polynomial of degree $r$ 
in the classical $\P r$-Lagrange basis on $\tref$.
The interpolant is denoted by $\ftr$ and we define
$\tr := \ftr(\tref)$.
The curved mesh of degree $r$ is~$\taur := \{ \tr; T \in \tauh \}$ 
with domain $\omhr := \cup_{\tr \in  \taur}\tr$ and with boundary $\ghr:= \partial \omhr$.
%
% and let $\phi_1^r, ... , \phi_{n_r}^r$ be the Lagrangian basis functions of degree $r \ge 1$ on~$\tref$ corresponding to the nodal points $\hat{x}^1, ... , \hat{x}^{n_r}$; where $n_r$ is the number of nodes on each element. The integer $ r \ge 1$ is referred to as the geometrical degree. Then we define the following diffeomorphism for each $T \in \tauh$ and each $r \ge 1$
% \begin{equation*}
% %\fonction{\ftr}{\tref  }{\tr:= \ftr(\tref)}{ \hat{x}}{ \displaystyle \ftr(\hat{x}) = \sum_{j=1}^{n_r} \fte(\hat{x}^j)\phi_j^r;}
% \ftr \, : \,  \hat{x} \in \tref  \longmapsto \sum_{j=1}^{n_r} \fte(\hat{x}^j)\phi_j^r \in \tr:= \ftr(\tref) ,
% \end{equation*}
% where $\fte:\tref \to \te$ is defined previously. Denote the higher order mesh $\taur := \{ \tr; T \in \tauh \}$ with $\omhr= \cup_{\tr \in  \taur}\tr$ as its domain and $\ghr:= \partial \omhr$ denotes its border. For more details, see \cite{elliott,ciaravtransf}. 
Note that $\ftr (v) = \fte (v)$ for 
$v$ a $\P r$-Lagrange node in $\tref$.
%any vertex  $v \in \tr$, by definition of $\ftr$.
\begin{exmp}
%\label{ex3}
%\emph{
In the quadratic case $r=2$ is displayed a border quadratic element $T^{(2)}$.
%In the quadratic case $r=2$ is displayed a border element $T$, the associated curved element $\te$ and the quadratic element $T^{(2)}$.
The mappings $F_T^{(2)}$ and $\fte$ coincide at the $\P 2$-Lagrange nodes which are the three vertexes $\hat{v}_i$ 
and the three edge mid-points $\hat{e}_i$ of $\tref$.
%
% Let $\tdeux$ be a quadratic border element ($\textcolor{red}{r = 2}$) constructed from a reference element $\tref$ with the help of the transformation $\ftdeux$, with $v^{(2)}_i = \ftdeux(\hat{v}_i)=~\fte (\hat{v}_i)$ and $e^{(2)}_ i~=~\ftdeux(\hat{e}_i)$ as the images of the vertices and the middle points of the edges respectively for $i=1, 2, 3$.
%}
\begin{center}
\vspace{-0.2cm}
\begin{tikzpicture}[scale=0.90]
%%%% conctruction of the reference element

\draw (0.5,0.5) node[above] {$\hat{T}$};

\draw (0,0) node  {$\bullet$};
\draw (2,0) node  {$\bullet$};
\draw (0,2) node  {$\bullet$};

\draw (0,0) node[below] {$\hat{v}_1$};
\draw (2,0) node[below]  {$\hat{v}_2$};
\draw (0,2) node[above]  {$\hat{v}_3$};

\draw (0,0) -- (2,0) ;
\draw (0,0) -- (0,2);
\draw (2,0) -- (0,2) ;

%%%%%%%%%%%%%%%%%%% curved ref element 

\draw[red] (1,0) node  {$\bullet$};
\draw[red] (1,1) node  {$\bullet$};
\draw[red] (0,1) node  {$\bullet$};

\draw (1,0) node[below] {$\hat{e}_3$};
\draw (1,1) node[above]  {$\hat{e}_1$};
\draw (0,1) node[left]  {$\hat{e}_2$};

%%%%%%%%%%%%%%%%%%%%%%%%%%%%%%%%%%%%% arrow and function

%\draw [->] (2.8,1) -- (4.3,1);
%\draw (3.5,1) node[above] {${\ft}$};

%%%%%%%%%%%%%%% curved case : transformation and red arrow

\draw[thick, red] [->] (2.8,1) -- (4.3,1);
\draw[red] (3.5,1) node[above] {$\textcolor{red}{\ftdeux}$};

%%%%%%%%%%%%%%%%%%%% triangle construction 

%%%%%%%%%%%%%%%%%%%%%%%%%%%%%%%%%%%%% name  

%%%%%%%%%%% linear
%\draw (7.5,0.8) node[above] {$T$};

%%%%%%%%%%% curved 
\draw[red] (7.4,0.8) node[above] {$\tdeux$};

%%%%%%%%%%%%%%%%%%%%%%%%%%%% image triangle

\draw (6,0) node  {$\bullet$};
\draw (7.36,2.5) node  {$\bullet$};
\draw (9,0) node  {$\bullet$}; 
\draw (7.36,2.56) -- (9,0);
\draw (9,0) -- (6,0);

%%%%%%%%%%%%% Linear 

%\draw (6,0)  node[below] {$v_1$};
%\draw (7.36,2.6) node[above]  {$v_3$};
%\draw (9,0) node[below]  {$v_2$};

%%%%%%%%%%%%%%%%%%%%%%%%%%% curved triangle

\draw (6,0)  node[left] {$v_1$};
\draw (7.36,2.6) node[above]  {$v_3$};
\draw (9,0) node[below]  {$v_2$};

\draw (6.3,1.3) node[left]  {$e_2$};
\draw (7.5,0) node[below]  {$e_3$};
\draw (8.15,1.3) node[right]  {$e_1$};

\draw[red] (7.5,0) node  {$\bullet$};
\draw[red] (6.4,1.5) node  {$\bullet$}; %(6.4,1.4) pour l'autre ex
\draw[red] (8.15,1.3) node  {$\bullet$};

%%%%%%%%%%%%%%%%%%%% (ARC edge) %%%%%%%%

\draw[red] (6,0) arc (180:122:3);

%%%%%%%%%%%%%%%%%%%%%%%%% the linear triangle

%\draw (6,0) -- (7.35,2.5);

%%%%%%%%%%%% x in T_1
%\draw (6.7,1.3) node  {$\bullet$};
%\draw (6.75,1.3) node[right]  {$x$};

%%%%%%%%%%% b(x)
%\draw (6.4,1.45) node  {$\bullet$};
%\draw (6.33,1.45) node[left]  {$b(x)$};

%%%%%%%%%% linking x and b(x)

%\draw[blue] (6.9,1.2) -- (6,1.7);

%%%%%%%%%%%%%%%%%%%%%% GAMMA NOT EXACT

\draw plot [domain=-0.2:1.5] (\x+6,-1.786*\x^2+4.3*\x);

\draw (5,-0.6) node {$\Gamma$};
\draw[thick]   [->] (5.2,-0.6)--(5.88,-0.5);

%%%%%%%%%%%%%%%%%%%% GAMMA (ARC) %%%%%%%%

%\draw (6,0) arc (180:90:3);
%\draw (9,3.2) node {$\Gamma$};

%%%%%%%%%%%%%%%%%%%%%% GAMMA_h

\draw[red] (5.7,2.75) node {$\Gamma_h^{(2)}$};
\draw[thick,red]   [->] (5.9,2.6)--(6.8,2);

\end{tikzpicture}
\end{center}
% \emph{Comparing $T$ to $\tdeux$ in the examples \ref{expl} and \ref{ex3}, one deduces that a quadratic element approximates $\Gamma$ better than a simplex, which is expected .}
\end{exmp}
\section{Numerical experiments}
\label{sec:num-exp}
%%
%%
%%
%%
%%
% In this section, we first validate the code Cumin developed by C. Pierre by considering the Laplace-Beltrami problem \eqref{lap} on a sphere. We estimate the error produced when approximating the solution of the equation by the finite element method. Afterwards, we consider the system \eqref{1} on the unit disk and calculate the norms of the error produced during the process.  
%%
%%
%%
%%
%%
%%
%%
\paragraph{Functional lift.}
Here, we define \textit{lifts} to transform a function on a domain $\omhr$ or $\ghr$ (defined in the previous section) into a function defined on $\Omega$ or $\Gamma$ respectively.
Lifts are necessary for two reasons: to compare the numerical solutions to the exact one and thus perform \textit{a priori} error estimates, but also to define the right hand side source terms in the numerical formulation of problem (\ref{1}). 

A surface lift is obviously provided by the orthogonal projection 
$b:~\ghr\rightarrow \Gamma$, 
to $v_h\in {\rm L}^2(\ghr)$ is associated $v_h^L \in {\rm L}^2(\Gamma)$  given by $v_h^L\circ b = v_h$.

To define a volume lift, a transformation 
$G_h^{(r)}:~ \omhr\rightarrow \Omega$
is defined and then to $u_h\in {\rm L}^2(\omhr)$ is associated $u_h^\ell \in {\rm L}^2(\Omega)$  given by $u_h^\ell\circ G_h^{(r)} = u_h$.
The definition of $G_h^{(r)}$ is less obvious and we describe it here.

In \cite{elliott}, it is given piecewise on all 
$\tr\in \taur$ by
${G_h}_{|_{\tr}} := \fte \circ ({\ftr})^{-1}$, where $T$ is the affine element relative to $\tr$.
However, this transformation does not fit the orthogonal projection~$b$ on the mesh boundary.
Precisely, following remark \ref{rem:Fe_T}, for $x\in\ghr\cap\tr$,
$G_h(x) := b \circ \ft\circ ({\ftr})^{-1}(x)$.
As a result the surface and bulk lifts do not coincide on 
$\ghr$: $  \left ( {\rm Tr} ~u_h\right )^L \neq {\rm Tr} (u_h^\ell)$.

To avoid this, 
we propose the following alternative definition of $G_h^{(r)}$ that is given piecewise for all 
$\tr\in \taur$ by
(with the notations of equation (\ref{eq:def-fte})),
\begin{equation}
  \label{eq:def-fter}
  {G_h^{(r)}}_{|_{\tr}} := \ftre \circ ({\ftr})^{-1},
  \quad 
  \ftre(\hat{x}) := x + (\lambda^\star)^{r+2} ( b(y) - y),
\end{equation}
with $x = \ftr(\hat{x})$ and $y = \ftr(\hat{y}).$
Geometrically, $\tr$ is directly transformed into $\te$ by $\ftre \circ (\ftr)^{(-1)}$, without being first transformed into $T$ as previously done.
Now, for $x\in \ftr\cap\ghr$, $\hat{x} = (\ftr)^{(-1)}(x)$ satisfies $\lambda^\star=1$ and so $\hat{y} = \hat{x}$ and $y=x$. 
So $\ftre(\hat{x}) = b(x) $, 
the volume and surface lifts both coincide with $b$ on $\ghr$ and the expected relation,
\begin{displaymath}
  \forall ~ u_h \in {\rm H}^1(\omhr), \quad 
  \left ( {\rm Tr} ~u_h\right )^L = {\rm Tr} (u_h^\ell),
\end{displaymath}
now holds.
Consequently, the surface lift $v_h^L$ now simply will be denoted by $v_h^\ell$.
\paragraph{Finite element formulation and implementation.}
On a mesh $\taur$
 is considered the finite element space, 
 \begin{equation}
   \label{eq:def-Vh}
   V_h := \left \{
u \in {\rm C}^0(\omhr), ~ \forall ~T\in \taur, 
~u_{\vert T} \circ \ftr \in \P k(\tref)
\right \} ,
 \end{equation}
with $\P k(\tref)$ the polynomials of degree $k$ on $\tref$ 
and with $k\ge 1$ the finite element degree.
Following \cite{elliott}, the problem \eqref{1} is discretized as: find $u_h\in V_h$ such that,
\begin{equation}\label{eq:fe-form}
  \forall  v_h\in V_h, \quad 
  a_h(u_h,v_h) = l_h(v_h) := \int_{\omhr}  f^{-\ell}\, J_{G_h^{(r)}}\d x
  + \int_{\ghr} g^{-\ell}\, J_b \d \sigma,
\end{equation}
with $G_h^{(r)}$ defined in (\ref{eq:def-fter}), 
with $f^{-\ell} := f\circ G_h^{(r)}$ 
and $g^{-\ell}:=g\circ b$ 
the inverse lifts of the source terms in \eqref{1}, 
with $J_{G_h^{(r)}}$ and  $J_b$ 
the Jacobians of $G_h^{(r)}$ and $b_{\ghr}$ respectively and where
$a_h$ is the bilinear form in (\ref{eq:def-a}) rewritten on $\omhr$ and $\ghr$.

Finite element space definition, matrix assembling and computation on curved surfaces are led 
%
% Entrée bibliographique à modifier.
%
using the code Cumin \cite{cumin}.
All integral computations rely on quadrature rules on the reference elements which are always chosen of sufficient order without further details.
\paragraph{Laplace equation on a surface.}
%\label{secdemlow}
%%
%%
%%
%%
In order to validate the code, we first 
draw our attention towards the Laplace equation 
$-\Delta_\Gamma u +u = g$
on a smooth surface $\Gamma\subset \R^3$.
We refer to Demlow~\cite{D1,D2} for the analysis of its finite element 
formulation.
Given a mesh $\taur$ of $\Gamma$, following \eqref{eq:def-Vh}, 
the~$\P k$-Lagrange finite element space is 
$W_h := 
\left \{
u \in {\rm C}^0(\ghr), ~ \forall ~T\in \taur, 
~u_{\vert T} \circ \ftr \in \P k(\tref)
\right \}$ and
the discrete problem is: find $u_h \in W_h$ such that,
\begin{displaymath}
  \forall ~ v_h \in W_h, 
  \quad 
  \int_{\ghr} \nabla_T u_h \cdot \nabla_T v_h \, \d \sigma
  +
  \int_{\ghr} u_h \, v_h \, \d \sigma = 
  \int_{\ghr} v_h \, g^{-\ell} \, J_b \,\d \sigma,
\end{displaymath}
with $g^{-\ell}$ and $J_b$ previously defined in (\ref{eq:fe-form}).
The \textit{a priori} error estimate for this problem developed by Demlow reads,
\begin{equation}
  \label{Laplace-surface-error-estimate}
    \|u-u_h^\ell \|_{\L2(\Gamma) } = O(h^{k+1} + h^{r+1}), \quad 
    \| \nt (u-u_h^\ell) \|_{\L2 (\Gamma) }  = O(h^k+h^{r+1}),
\end{equation}
for a smooth enough source term $g$.

We set $\Gamma$ to the unit sphere 
and the source term to $g(x,y,z) = \mathrm{e}^y (y+2)y$.
Three series of successively refined meshes, respectively affine, quadratic and cubic, of $\Gamma$ have been generated by the software 
Gmsh\footnote{Gmsh: a three-dimensional finite element mesh generator, \url{https://gmsh.info/}}.
The numerical errors have been computed 
for each mesh  and 
for $\mathbb{P}^k$, with $k=1,\dots 4$.
\begin{figure}[h!]
\centering
    \includegraphics[width=0.3\textwidth]{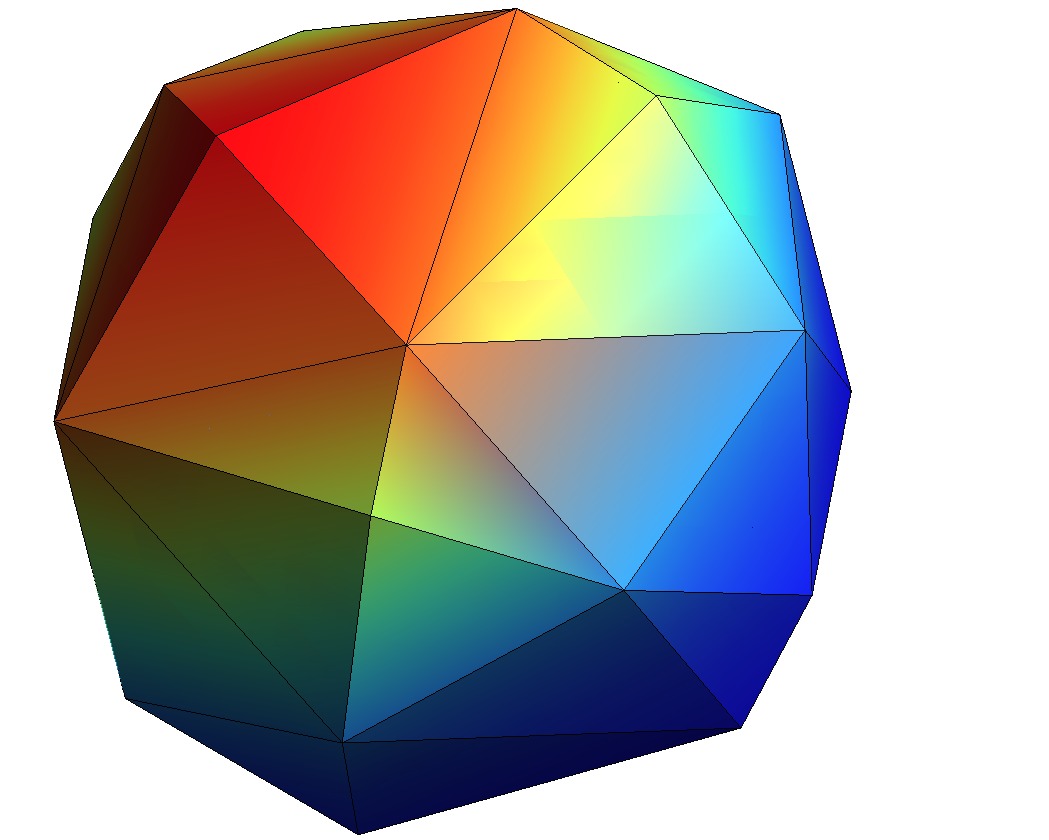}
    ~~~
    \includegraphics[width=0.3\textwidth]{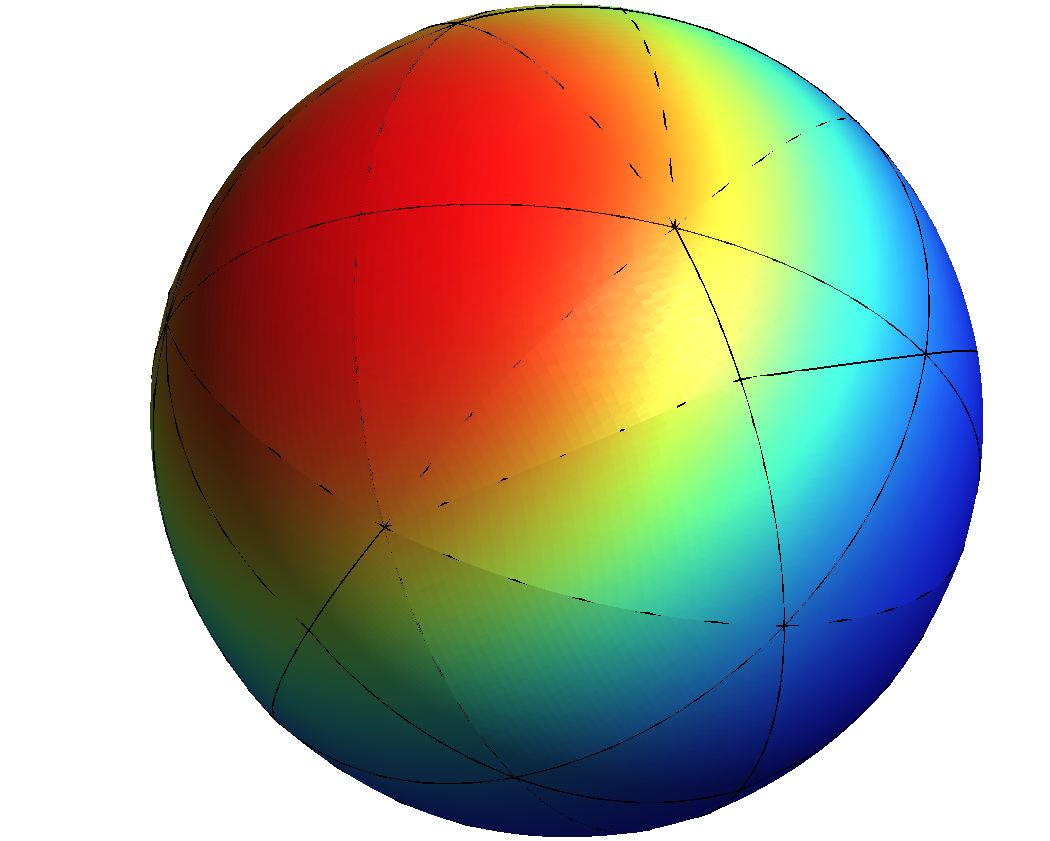}
 \caption{Numerical solution of 
the Laplace equation on a sphere with affine and quadratic meshes.}
 \label{f1}
\end{figure}
\begin{table}[!ht]
    \centering
    \begin{tabular}{|l||l|l|l|l||l|l|l|l|}%{|l|M|C|R|J|}
\cline{2-9}
\multicolumn{1}{c||}{}    &  \multicolumn{4}{|c||}{$\| u-u_h^\ell \|_{\L2 (\Gamma)}$} & \multicolumn{4}{|c|}{$\| u-u_h^\ell \|_{\HH1 (\Gamma)}$}  \\[0.08cm]
\cline{2-9}
\multicolumn{1}{c||}{}    &  $\P1$ &  $\P2$ &  $\P3$ &  $\P4$ & $\P1$ &  $\P2$ &  $\P3$ &  $\P4$  \tabularnewline
\hline
 Affine mesh (r=1)   & \textcolor{black}{1.96} & \textcolor{black}{1.96} & \textcolor{black}{1.96} & \textcolor{black}{1.96} & 0.99 & \textcolor{black}{1.96} & \textcolor{black}{1.96} & \textcolor{black}{1.96} \tabularnewline
\hline
Quadratic mesh (r=2) & 1.98 & 2.95 & \textcolor{red}{3.92} & \textcolor{red}{3.92} & 0.98 & 1.97 & 3.00 & \textcolor{red}{3.91} \tabularnewline
\hline
Cubic mesh (r=3) & 1.98 & 2.94 & \textcolor{black}{3.95} & \textcolor{black}{3.92} & 0.98 & 1.96 & 2.96 & \textcolor{black}{3.95} \tabularnewline
\hline
\end{tabular}
\caption{\label{TABLE1}Convergence order for the Laplace equation on a sphere.}
\end{table}

The numerical solution on two coarse meshes is depicted on figure \ref{f1},
and the measured convergence orders are reported in table \ref{TABLE1}.
The affine and cubic meshes behave exactly as expected following 
\eqref{Laplace-surface-error-estimate}.
In turn, quadratic meshes produce unexpected convergence rates indicated in red in table \ref{TABLE1} and a super convergence is observed.
Quadratic meshes display a geometrical error $h^4$ instead of the expected $h^3$ and thus behave as if $r=3$.
This behavior has been further investigated and is not problem dependent. It is also observed for the Poisson problem on a disk with Neumann or Robin boundary conditions. 
It is neither caused by the considered geometry: 
studying a simpler problem of integral computation on a non-symmetric and non-convex domain gave the same surprising super convergence.
So far we have no further explanation for this particular error.
\paragraph{Numerical study of the Ventcel problem.} 
The Ventcel problem (\ref{1}) is considered on the unit disk $\Omega$
with $\alpha=\beta=1$ and~$\kappa=0$, 
with the source terms $f(x,y)= - y \mathrm{e}^{x}$ and 
$g(x,y)=  y \mathrm{e}^{x} ( 3 + 4x-y^2 )$
corresponding to the exact solution
$u = -f$.
The discrete problem (\ref{eq:fe-form}) 
is implemented 
%
% Entrée bibliographique à modifier.
%
and solved using the code Cumin \cite{cumin}.
Again, three series of successively refined meshes, respectively affine, quadratic and cubic, of $\Omega$ have been generated with Gmsh.
For each mesh and for  $\mathbb{P}^k$ finite elements, with $k=1\dots 4$,
four numerical errors are computed (two in the bulk domain and two on the boundary),
\begin{displaymath}
  \|u-u_h^\ell \|_{\L2(\Omega) }, \quad 
  \| \na (u-u_h^\ell) \|_{\L2 (\Omega) },
  \quad 
  \|u-u_h^\ell \|_{\L2(\Gamma) } \quad {\rm and} \quad 
  \| \nt (u-u_h^\ell) \|_{\L2 (\Gamma) },
\end{displaymath}
and the estimated convergence rates are reported in the two tables 
\ref{TABLE3} and \ref{TABLE2}.
\begin{table}[!ht]
    \centering
    \begin{tabular}{|l||l|l|l|l||l|l|l|l|}%{|l|M|C|R|J|}
\cline{2-9}
\multicolumn{1}{c||}{}    &  \multicolumn{4}{|c||}{$\| u-u_h^\ell \|_{\L2 (\Gamma)}$} & \multicolumn{4}{|c|}{$\| \nabla _\Gamma (u- u_h^\ell) \|_{\L2 (\Gamma)}$}  \\[0.08cm]
\cline{2-9}
\multicolumn{1}{c||}{}    &  $\P1$ &  $\P2$ &  $\P3$ &  $\P4$ & $\P1$ &  $\P2$ &  $\P3$ &  $\P4$  \tabularnewline
\hline
 Affine mesh (r=1)   & 
\textcolor{black}{2.00} & \textcolor{black}{2.03} & \textcolor{black}{2.01} & \textcolor{black}{2.01} & 
1.00 & \textcolor{black}{2.00} & \textcolor{black}{1.99} & \textcolor{black}{1.98} \tabularnewline
\hline
Quadratic mesh (r=2) & 
2.00 & 3.00 & \textcolor{black}{4.00} & \textcolor{black}{4.02} & 
1.00 & 2.00 & 3.00 & \textcolor{black}{4.02} \tabularnewline
\hline
Cubic mesh (r=3) & 
2.00 & 3.00 & 4.00 & 4.24 & 
1.00 & 2.00 & 3.00 & 3.98 \tabularnewline
\hline
\end{tabular}
\caption{\label{TABLE3}Convergence order of $\| u-u_h^\ell \|_{\L2 (\Gamma) } $ and of 
$\| \nabla _\Gamma (u- u_h^\ell) \|_{\L2 (\Gamma)}$}
\end{table}
The surface errors in table \ref{TABLE3} behave exactly the same way 
as the estimation
\eqref{Laplace-surface-error-estimate} for 
the Laplace equation on a surface: %, see equation:
the same super-convergence for the quadratic meshes again occurs, as if $r=3$ in that case.
As a consequence, the numerical solution seems to be correctly computed.

\begin{table}[!ht]
    \centering
    \begin{tabular}{|l||l|l|l|l||l|l|l|l|}%{|l|M|C|R|J|}
\cline{2-9}
\multicolumn{1}{c||}{}    &  \multicolumn{4}{|c||}{$\| u-u_h^\ell \|_{\L2 (\Omega)}$} & \multicolumn{4}{|c|}{$\| \nabla  (u-  u_h^\ell) \|_{\L2 (\Omega)}$}  \\[0.08cm]
\cline{2-9}
\multicolumn{1}{c||}{}    &  $\P1$ &  $\P2$ &  $\P3$ &  $\P4$ & $\P1$ &  $\P2$ &  $\P3$ &  $\P4$  \tabularnewline
\hline
 Affine mesh (r=1)   & 
\textcolor{black}{1.98} & \textcolor{black}{1.99} & \textcolor{black}{1.97} & \textcolor{black}{1.97} & 
1.00 & \textcolor{black}{1.50} & \textcolor{black}{1.49} & \textcolor{black}{1.49} \tabularnewline
\hline
Quadratic mesh (r=2) & 
2.01 & 3.14 & \textcolor{black}{3.94} & \textcolor{black}{3.97} & 1.00 & 2.12 & 3.03 & \textcolor{black}{3.48} \tabularnewline
\hline
Cubic mesh (r=3) & 
2.04 & \textcolor{red}{2.45} & \textcolor{red}{3.44} & 4.04 & 
1.02 & \textcolor{red}{1.47} & \textcolor{red}{2.42} & 3.46 \tabularnewline
\hline
\end{tabular}
\caption{\label{TABLE2}Convergence order of $\| u-u_h^\ell \|_{\L2 (\Omega) } $ and $\| \nabla (u- u_h^\ell) \|_{\L2 (\Omega) } $}
\end{table}
The interpretation of the convergence rates for the bulk errors 
in table \ref{TABLE2} is less straightforward.
Let us first focus on the affine and quadratic meshes, and consider that in the quadratic case $r=3$ instead of $2$ (as a consequence of the super convergence in that case previously discussed).
Then the figures in table \ref{TABLE2} can be interpreted as,
\begin{displaymath}
  \| u-u_h^\ell \|_{\L2 (\Omega) } = O( h^{k+1} + h^{r+1})
\quad  {\rm and } \quad  
  \| \nabla (u- u_h^\ell) \|_{\L2 (\Omega) } = O( h^{k} + h^{r+1/2}).
\end{displaymath}
This behavior differs from (\ref{Laplace-surface-error-estimate}) for the gradient norm where $h^{r+1}$ is now replaced by $h^{r+1/2}$.
This difference could be understood from a theoretical point of view 
following ideas that should be presented in fore coming works.

For the cubic case, the $\P1$ and $\P4$ cases behave the same way.
However, for the $\P2$ and~$\P3$ cases (red figures in table \ref{TABLE2}), 
the rule rather seem to be $\| u-u_h^\ell \|_{\L2 (\Omega) } = O( h^{k+1/2} + h^{r+1})$ and~$\| \nabla u-\nabla u_h^\ell \|_{\L2 (\Omega) } = O( h^{k-1/2} + h^{r+1/2})$. 
Though we have no clear understanding on this, we experienced that the choice of the lift operator here has  a crucial influence.
We recall that this lift is based on a geometric transformation 
$G_h^{(r)}:~ \omhr\rightarrow \Omega$, which is a modification of the one
defined in \cite{elliott}.
%\red{We recall that this lift is based on a geometric transformation $G_h^{(r)}:~ \omhr\rightarrow \Omega$ that we adopted a modification of that transformation, defined in \cite{elliott}.}
When resorting to the lift in \cite{elliott}, a saturation of the convergence order is observed: 2.5 for the $\L2 (\Omega)$ norm  and 1.5 for the gradient $\L2$-norm on $\Omega$.
The same observation holds both for the quadratic and cubic meshes.

% As before, we observe an error saturation at $\textcolor{blue}{r+1=2}$ in the affine case (r=1), for the surface norms and the $\L2(\Omega)$ norm of the error. However, we see that $\| u-u_h^\ell \|_{\HH1(\Omega)}$ saturates at~$\textcolor{blue}{h^{r+1/2}=h^{1.5}}$. On the quadratic mesh, the error norms all display \textcolor{red}{better results than expected} (in red in the tables), as seen in the previous example. As before, these results suggest that, for~$r=2$, the error behaves as if it is approximated on a cubic mesh (r=3). Lastly, for $r=3$, the surface norms of the error saturate as expected at $\textcolor{blue}{r+1=4}$. However, we observe unexpected results for the volume norms. For both norms, the convergence order of the error seem to saturate at $\textcolor{green}{k-1/2}$ (in green in the tables) which is not a common result. 

% \medskip

% \begin{rem}
% It has to be noted that there are some numerical difficulties, especially when constructing the mesh. The definition of the transformation $G_h$ that appears in the lift operator above and its derivatives is not simple to implement. To be able to estimate the error between the exact solution $u$ and its finite element approximation $u_h$ defined on separate domains $\Omega$ and $\omhr$ respectively, the lift is always required to be well defined. Thus, the transformation~$G_h$ needs to be correctly defined.
% \end{rem}
%%
%%
%%
%%
%%
%%
%%
%%
\section*{Conclusion}
%We presented a theoretical study of the Ventcel problem 
%(\ref{1}) together with a numerical analysis of the associated 
%\textit{a priori} errors using high order finite elements on curved meshes.
%This numerical analysis is supported by an alternative definition of a lift operator as compared to \cite{elliott} which improved our numerical results.
%{\color{red}
We have presented an approach in order to numerically solve 
the Ventcel problem~\eqref{1} and have used the code Cumin~\cite{cumin} 
to give a numerical exploration of the associated 
\textit{a priori} errors using high order finite elements on curved meshes.
This numerical analysis is supported by an alternative definition of a lift operator 
as compared to the previous work~\cite{elliott} which improved our numerical results.
Beyond difficulties related to the lift definition, and beyond unexplained super convergence associated to quadratic meshes, we formulate the following conjecture for the Ventcel problem \textit{a priori} numerical errors, 
the proof of which is a work in progress.
\\[5pt]
\textbf{Conjecture.}
Let $u \in \Hk1\omgam$ be a solution of the variational problem \eqref{fv1}, 
let $\taur$ be a mesh of $\Omega$ with geometrical degree $r$,
let $V_h$, defined in \eqref{eq:def-Vh}, be the associated finite 
element space of degree $k$. Then the numerical solution $u_h\in V_h$ to
the discrete problem (\ref{eq:fe-form}) satisfies,
\begin{align*}
  &\| u-u_h^\ell \|_{\L2 (\Omega) } = O( h^{k+1} + h^{r+1}),
  \quad  \quad  \quad  ~
  \| \nabla (u- u_h^\ell) \|_{\L2 (\Omega) } = O( h^{k} + h^{r+1/2}),
  \\
    &\|u-u_h^\ell \|_{\L2(\Gamma) } = O(h^{k+1} + h^{r+1}), \quad 
  {\rm and } \quad  
    \| \nt (u-u_h^\ell) \|_{\L2 (\Gamma) }  = O(h^k+h^{r+1}).
\end{align*}

\bibliographystyle{abbrv}
\bibliography{biblio}

\end{document}